\documentclass[aos,MSNbibl,nameyear,dvips]{arximspdf}
\usepackage{dcolumn}
\usepackage{graphicx}

\doi{10.1214/13-AOS1129} 
\volume{41}
\issue{3}
\pubyear{2013}
\firstpage{1669}
\lastpage{1692}

\makeatletter

\newcolumntype{d}[1]{D{.}{.}{#1}}

\newtheorem{theorem}{Theorem}[section]
\newtheorem{lemma}{Lemma}[section]

\newcommand{\sE}{\mathsf{E}}
\newcommand{\var}{\operatorname{var}}
\newcommand{\mse}{\operatorname{mse}}
\newcommand{\cov}{\operatorname{Cov}}

\newcommand{\btheta}{\bolds{\theta}}

\newcommand{\diag}{\operatorname{diag}}
\newcommand{\bq}{\mathbf{q}}
\newcommand{\bh}{\mathbf{h}}
\newcommand{\bB}{\mathbf{B}}

\newcommand{\bH}{\mathbf{H}}
\newcommand{\bW}{\mathbf{W}}
\newcommand{\bZ}{\mathbf{Z}}

\newcommand{\bS}{\mathbf{S}}
\newcommand{\bI}{\mathbf{I}}
\newcommand{\0}{\mathbf{0}}

\newcommand{\EL}{\operatorname{EL}}
\newcommand{\EM}{\operatorname{EM}}
\makeatother

\begin{document}
\begin{frontmatter}

\title{Quantile and quantile-function estimations under density ratio
model\thanksref{T1}}
\thankstext{T1}{Supported in part by a Collaborative Research and Development
grant from Natural Science and Engineering Research Council of Canada
and in-kind contributions from Natural Resources Canada (Canadian
Forest Service)
through FPInnovation's national research program.}
\runtitle{Quantile and quantile-function estimations}

\begin{aug}
\author[A]{\fnms{Jiahua} \snm{Chen}\ead[label=e1]{jhchen@stat.ubc.ca}}
\and
\author[B]{\fnms{Yukun} \snm{Liu}\corref{}\ead[label=e2]{ykliu@sfs.ecnu.edu.cn}\thanksref{t2}}
\thankstext{t2}{Supported by the NNSF of China,
11001083 and 11101156.}
\runauthor{J. Chen and Y. Liu}
\affiliation{University of British Columbia and East China Normal University}
\address[A]{Department of Statistics\\
University of British Columbia\\
Vancouver, BC\\
Canada V6T 1Z2\\
\printead{e1}}

\address[B]{Department of Statistics and Actuarial Sciences\\
School of Finance and Statistics\\
East China Normal University\\
Shanghai 200241\\
China\\
\printead{e2}}
\end{aug}

\received{\smonth{8} \syear{2012}}
\revised{\smonth{4} \syear{2013}}


\begin{abstract}
Population quantiles and their functions are important
parameters in many applications. For example,
the lower quantiles often serve as crucial quality
indices for forestry products.
Given several independent samples from
populations satisfying the density ratio model,
we investigate the properties of empirical likelihood (EL) based inferences.
The induced EL quantile estimators are shown to
admit a Bahadur representation that leads to
asymptotically valid confidence intervals
for functions of quantiles. We rigorously prove
that EL quantiles based on all the samples are
more efficient than empirical quantiles
based on individual
samples. A simulation study shows that the EL quantiles and
their functions have superior performance when the density
ratio model assumption is satisfied and when it is mildly violated.
An example is used to demonstrate the
new method and the potential cost savings.
\end{abstract}

%
\begin{keyword}[class=AMS]
\kwd[Primary ]{62G20}
\kwd[; secondary ]{62G15}
\end{keyword}

\begin{keyword}
\kwd{Asympotic efficiency}
\kwd{Bahadur representation}
\kwd{confidence interval}
\kwd{empirical likelihood}
\end{keyword}

\end{frontmatter}

\section{Introduction}\label{sec1}
Forestry plays a major role in the Canadian economy;
maintaining the high quality of wood products is vital economically
and socially. We are designing an effective long-term monitoring
plan for the quality of forestry products in Canada.
Two important quality indices for a piece of lumber
are the modulus of elasticity (MOE) and the modulus of rupture (MOR),
its strength in terms of elasticity and toughness.
The reliability of lumber-based structures may depend heavily
on the lower population quantiles of these indices.
However, it is costly, time consuming and laborious
to obtain these quality measurements.
Therefore, efficient estimates of the population quantiles
and their functions are important.

The estimation of quantiles based on a single random sample
is a well-researched topic.
Empirical quantiles have been shown to admit a Bahadur representation
[\citet{Bah66}; \citet{Kie67}; \citet{Ser80}],
making it simple to study the joint limiting
distributions of any number of sample quantiles and their smooth functions.
In the presence of auxiliary information,
the empirical likelihood [EL; Owen (\citeyear{Owe88}, \citeyear{Owe01})]
can be utilized to improve efficiency.
The Bahadur representation of EL
estimators has been established by \citet{CheChe00}.
There is also an abundant literature on the Bahadur representation
when the samples are not independent or
have a time-series structure. See \citet{Wu05}
and \citet{ZhoWu09} for recent examples.


In the targeted application,
we potentially have a number of random samples from
similar populations, and the combined sample size is large.
Even if the size of each random sample is small,
the total sample size increases over time.
We may also have samples from similar products, such as lumber of
various shapes and lengths.
If the population distributions have some common features,
the pooled information may improve the efficiency
of each quantile estimate.

Specifically, we study quantile
estimators based on the density ratio model (DRM) of \citet{And79}
and the EL approach, and we focus on investigating their properties.
Suppose we have $m+1$ independent random samples from populations
with cumulative distribution and density functions denoted
$G_k(x)$ and $g_k(x)$, $k = 0, 1, \ldots, m$.
The DRM postulates that
%
\begin{equation}
\label{DRM} \log\bigl\{g_k(x)/g_0(x)\bigr\} =
\btheta_k^{\tau} \bq(x)
\end{equation}
for some known function $\bq(x)$ of dimension $d$ and
corresponding unknown vector-valued parameters $\btheta_k$.
We require the first element of $\bq(x)$ to be one so that
the first element of $\btheta_k$ is a normalization parameter.

In this formulation, the form of $g_0(x)$ is unspecified.
Many parametric distribution families including normal and Gamma
are special cases of the DRM.
\citet{QinZha97} showed that the
logistic regression model commonly used in case--control studies can
be described by the DRM. They studied the EL approach for
parameter estimation and for goodness-of-fit tests of the regression model.
Zhang (\citeyear{Zha00,Zha02}) investigated the EL approach
for quantile estimation
and goodness-of-fit.
\citet{Foketal01} used the EL approach under the DRM
for a classical one-way analysis-of-variance.

We focus on DRM-based quantile estimation
and study its Bahadur representation. We show
that the EL quantiles are more efficient than empirical quantiles.
The representation is then used to construct confidence intervals
for the quantiles and their functions.
These results are particularly relevant for the design of a long-term
monitoring system
for wood products. The finite-sample performance of the new methods
is superior to that of the empirical quantiles when
the DRM model assumption is valid and when it is mildly violated.

In Section~\ref{sec2}, we review the EL approach under the DRM.
Section~\ref{sec3} derives the Bahadur representation.
In Section~\ref{sec4}, we study the asymptotic properties
of the new quantile estimates.
The finite-sample performance is examined in Section~\ref{secquancom},
and the proposed quantile estimation is illustrated using
lumber data in Section~\ref{sec6}.
The proofs are given in the \hyperref[app]{Appendix}.

\section{Empirical likelihood under DRM}\label{sec2}

Empirical likelihood under the DRM
can be found in \citet{QinZha97}
or \citet{Foketal01}.
Suppose the population distributions, $G_k(x)$, of
$m + 1$ random samples of sizes $n_k$:
$
\{ (x_{kj} \dvtx j = 1, \ldots, n_k); k = 0, 1, 2,
\ldots, m \}
$
satisfy the DRM (\ref{DRM}).
The model assumption may also be written
\[
dG_k(x) = \exp\bigl\{ \btheta^{\tau}_k \bq(x)
\bigr\} \,dG_0(x).
\]
If $G$ is discrete, then
$dG(x) = G(x) - G(x_-) = P(X = x)$ for the corresponding
random variable $X$. The EL
is defined as if the $G$'s are discrete.
Let $p_{kj} = dG_0(x_{kj})$ for all $k, j$.
The EL is defined as
\[
L_n(G_0, G_1, \ldots, G_m) =
\prod_{k, j} \,dG_k(x_{kj}) =
\biggl\{ \prod_{k, j} p_{kj} \biggr\} \times
\exp\biggl\{ \sum_{k, j} \btheta_k^{\tau}
\bq(x_{kj}) \biggr\},
\]
where the product and summation with respect to
$\{k, j\}$ are over the full range:
$k=0, \ldots, m$ and $j=1, \ldots, n_k$.
We set $\btheta_0 = 0$ for notational simplicity.

The DRM assumption implies that $L_n$ is
also a function of the parameter vector
$\btheta^\tau= (\btheta_1^\tau, \ldots, \btheta_m^\tau)$ and $G_0$.
Hence, we may also write its logarithm as
\[
\ell_n (\btheta, G_0) = \sum
_{k, j} \log(p_{kj}) + \sum
_{k, j} \btheta_k^{\tau}
\bq(x_{kj}).
\]
The model assumption also implies that, for $r =0, 1, \ldots, m$,
%
\begin{equation}
\label{constr} \int\exp\bigl\{\btheta_r^{\tau} \bq(x)\bigr
\}\,dG_{0}(x) = 1.
\end{equation}
Thus, for any $r$ between 0 and $m$,
$
\sum_{k, j} p_{kj} \exp\{
\btheta_r^{\tau} \bq(x_{kj}) \} = 1,
$
which is naturally accommodated in the EL approach.

Inference on $\btheta$ and other aspects of the population distributions
is usually carried out by first profiling the EL with respect to
$\btheta$.
That is, we define
$\tilde\ell_n(\btheta) = \max_{G_0} \ell_n(\btheta, G_0)$
subject to constraints (\ref{constr}) on $G_0$.
Technically, we confine the support of $G_0$
to $\{ x_{kj}\}$.
The maximum in $G_0$ is attained when
\[
p_{kj} = n^{-1} \Biggl\{1+ \sum_{s=1}^m
\nu_s \bigl[\exp\bigl\{\btheta_s^{\tau}
\bq(x_{kj}) \bigr\}- 1 \bigr] \Biggr\}^{-1},
\]
where $(\nu_1, \nu_2, \ldots, \nu_m)$ is the solution to
$
\sum_{k, j} p_{kj} \exp\{
\btheta^\tau_r \bq(x_{kj}) \} = 1
$
for $r=1, 2, \ldots, m$
and $n= \sum_{k=0}^m n_k$ is the total sample size.
The profile log-EL
(up to an additive constant)
is given by
%
\begin{equation}
\label{truelik} \tilde\ell_n (\btheta) = - \sum
_{k, j} \log\Biggl\{ 1+ \sum_{s=1}^m
\nu_s \bigl[\exp\bigl\{\btheta_s^{\tau}
\bq(x_{kj})\bigr\}- 1\bigr] \Biggr\} + \sum
_{k, j} \btheta_k^{\tau}
\bq(x_{kj}).
\end{equation}

We may regard $\tilde\ell_n (\btheta)$ as
a parametric likelihood and activate the classical
likelihood-based statistical inference.
This profile likelihood has the
same maximum value and point as another
function,
%
\begin{equation}
\label{proflik} \ell_n(\btheta) = - \sum
_{k, j} \log\Biggl[\sum_{r=0}^m
\rho_{r} \exp\bigl\{\btheta_r^{\tau}
\bq(x_{kj})\bigr\} \Biggr] + \sum_{k, j}
\btheta_k^{\tau} \bq(x_{kj}),
\end{equation}
with $\rho_{r} = n_r/n$.
Because of its simplicity,
the literature often regards $\ell_n(\btheta)$ instead of
$\tilde\ell_n(\btheta)$
as the profile likelihood function of $\btheta$.
In the two-sample situation, \citet{KezLeo08}
found that it is a ``dual likelihood.''
Its likelihood ratio statistics remain asymptotically chi-square.

We define the maximum EL estimator (MELE)
$\hat{\btheta}$ of $\btheta$ as the maximum point of~(\ref{truelik})
or equivalently of (\ref{proflik}).
The asymptotic normality of $\hat{\btheta}$
has been established in various situations
[\citet{QinZha97}; \citet{Foketal01}; \citet{Zha02}].
We summarize and extend these results, giving the necessary details as a
preparational step. Let
$
h(x; \btheta) = \sum_{k=0}^m
\rho_{k} \exp\{ \btheta_k^{\tau} \bq(x)
\},
$
and for $k=0, 1, \ldots, m$, we define
\[
h_k(x; \btheta) = \rho_{k} \exp\bigl\{
\btheta_k^{\tau} \bq(x) \bigr\}/h(x; \btheta).
\]
Let
$\bh(x; \btheta) = \{ h_0(x; \btheta), \ldots, h_m(x; \btheta)\}^\tau$
and define an $(m+1) \times(m+1)$ matrix
\[
\bH(x; \btheta) = \diag\bigl\{ \bh(x; \btheta) \bigr\} - \bh(x;
\btheta)
\bh^\tau(x; \btheta).
\]
When $\btheta= \btheta^*$, the true value of $\btheta$,
we may drop $\btheta^*$
for notational simplicity.
Finally, we use $d\bar G(x)$ for $h(x; \btheta^*) \,dG_0(x)$ in the
integrations.

%
\begin{theorem}
\label{thetaest}
Suppose we have an independent random sample
$\{ x_{kj} \}_{j=1}^{n_k}$ from population $G_k$
for $k=0, 1, \ldots, m$.
The total sample size $n = \sum_k n_k \to\infty$,
and
$\rho_k = n_k/n$ remains a constant (or within the $n^{-1}$ range).

The population distributions $G_k$ satisfy the DRM (\ref{DRM}) with true
parameter value $\btheta^*$ and
$\int h(x; \btheta) \,d G_0 < \infty$ in a neighborhood of $\btheta^*$.
The components of $\bq(x)$ are linearly independent
and its first element is one.

Then $\sqrt{n} ( \hat{\btheta} - \btheta^*)$ is asymptotically
multivariate normal with mean $\0$
and covariance matrix $\bW^{-1} - \bS$.
Both $\bW$ and $\bS$ are $md \times md$
block matrices with each block a $d\times d$ matrix,
and their $(r, s)$th blocks are, respectively,
\begin{eqnarray*}
\bW_{rs} &=& \int\bq(x)\bq^{\tau}(x) \bigl\{ h_r
(x) \delta_{rs} - h_r(x)h_s(x) \bigr\} \,d\bar
G(x),
\\
\bS_{rs} &=& \bigl(\rho_r^{-1}
\delta_{rs}+\rho_0^{-1}\bigr) \operatorname{diag}
\{1, 0, \ldots, 0\},
\end{eqnarray*}
where $1\leq r, s \leq m$ and $\delta_{rs}=1$ if $r=s$ and 0 otherwise.
\end{theorem}

The \hyperref[app]{Appendix} contains a sketched proof to bridge some gaps.
Using the Kronecker product $\otimes$, we have a
tighter expression,
%
\begin{equation}
\label{Wdef} \bW= \int\bH{[-1, -1]}(x) \otimes\bigl\{
\bq(x)
\bq^\tau(x) \bigr\} \,d\bar G(x),
\end{equation}
where
$\bH{[-1, -1]}(x)$ is $\bH(x)$
with its first row and column removed.
This convention is adopted from the statistical software package
[\citet{Tea11}].

%

The asymptotic normality is a stepping stone for our main result on the
Bahadur representation. It also reveals that the MELE is root-$n$ consistent.
The assumption that
$\int h(x; \theta) \,dG_0 < \infty$ in a neighborhood
of $\btheta^*$ implies the existence of the
moment generating function of $\bq(x)$
and therefore all its finite moments.
This fact will be used in our proofs.

\section{Bahadur representation and its applications}\label{sec3}

Given the MELE $\hat{\btheta}$, the fitted values of $p_{kj}$ are
$
\hat p_{kj} = \{ n h(x_{kj}; \hat{\btheta})
\}^{-1}
$
and the fitted $G_r$ is
\begin{eqnarray*}
\hat G_r(x) &=& \sum_{k, j} \hat
p_{kj} \exp\bigl\{ \hat{\btheta}_r^\tau
\bq(x_{kj}) \bigr\} I(x_{k j} \leq x) \\
&=& n_r^{-1}
\sum_{k, j} h_r(x_{kj}; \hat{
\btheta}) I(x_{k j} \leq x)
\end{eqnarray*}
with $\hat{\btheta}_0 = 0$ and $I(A)$ an indicator function of event $A$.
For any $ \alpha\in(0, 1)$, define the $\alpha$-quantile of $G_r$ as
$\xi_r = \xi_{r,\alpha} = \inf\{ x\dvtx G_r(x)\geq\alpha\}$
and its EL-based estimator as
%
\begin{equation}
\label{ELquantile} \hat\xi_r = \hat\xi_{r, \alpha} = \inf\bigl
\{ x\dvtx\hat G_r(x)\geq\alpha\bigr\}.
\end{equation}
We call the $\hat\xi_r$ values EL quantiles for simplicity.
The asymptotic normality of the EL quantile $\hat\xi_{r}$
is useful for constructing confidence intervals for $\xi_{r}$ or
for testing related hypotheses.
Researchers are often interested in smooth functions of
quantiles of many populations and/or at several levels.
Thus, the multivariate asymptotic behavior is useful,
and this calls for the Bahadur representation.

%
\begin{theorem}
\label{bahadur}
Assume the conditions of Theorem~\ref{thetaest},
and that the density function $g_r(x)$ is continuous and
positive at $x = \xi_r$.
The EL quantile (\ref{ELquantile}) has Bahadur representation
%
\begin{equation}
\label{bhrep} \hat\xi_{r} = \xi_{r} + \bigl\{{\alpha} -
\hat G_r(\xi_{r})\bigr\}/g_r(
\xi_{r}) + O_p\bigl(n^{-3/4} \bigl\{\log(n)\bigr
\}^{1/2}\bigr),
\end{equation}
where $\xi_r$ is the {$\alpha$th quantile} of $G_r(x)$.
\end{theorem}

The proof is given in the \hyperref[app]{Appendix}.
Without the Bahadur representation, it is a daunting
task to derive the limiting distribution of
functions of EL quantiles such as
$\hat\xi_{0, \alpha} - \hat\xi_{1, \alpha}$.
Theorem~\ref{bahadur} links this task to that of
$\{a \hat G_0(\xi_{0, \alpha}) - b \hat G_1(\xi_{1, \alpha})\}$
with nonrandom constants $a$ and $b$.
The asymptotic properties of
$\hat G_r$ are simple and easy to use.

%
\begin{theorem}
\label{asyofdistest}
Assume the same conditions as in Theorem~\ref{thetaest}.
For any $0 \leq r_1, r_2, \ldots, r_k \leq m$ and an accompanying
set of real numbers $x_1, x_2, \ldots, x_k$ in the support of $G_0(x)$,
$\sqrt{n}\{\hat G_{r_j} (x_j) - G_{r_j}(x_j)\}$ are
jointly asymptotically $k$-variate normal with mean $\0$
and covariance matrix
$
\Omega_{\EL} = ( \omega_{r_i, r_j} ( x_i,
x_j ) )_{ 1 \leq i, j \leq k}
$.
The generic form of $ \omega_{r_i, r_j} ( x_i, x_j )$ is given by
%
\begin{equation}
\label{omegars} \omega_{rs}(x,y) = \sigma_{rs}(x, y) - (
\rho_r\rho_s)^{-1} \bigl\{ a_{rs}(x
\wedge y) - \bB_r^{\tau}(x)\bW^{-1}
\bB_s(y) \bigr\},
\end{equation}
where $x \wedge y = \min\{x,y\}$,
%
\begin{eqnarray*}
\sigma_{rs}(x, y) &=& \rho_r^{-1}
\delta_{rs} \bigl\{ G_r(x\wedge y) - G_r(x)G_s(y)
\bigr\},
\\
a_{rs}(x) &=& \int_{-\infty}^x \bigl\{
\delta_{rs} h_r(t) - h_r(t)
h_s(t) \bigr\} \,d\bar G(t)
\end{eqnarray*}
and $\bB_r(x)$ is a vector of length $md$ with its
$s$th segment (of length $d$ with $s=1, 2, \ldots, m$) being
\[
\bB_{r,s}(x) = \int_{-\infty}^x \bigl\{
\delta_{rs} h_r(t) - h_r(t)
h_s(t) \bigr\} \bq(t) \,d\bar G(t).
\]
\end{theorem}

The proof is given in the \hyperref[app]{Appendix}.
The Bahadur representation (\ref{bhrep})
and the multivariate asymptotic normality of the $\hat G$'s
lead to multivariate asymptotic normality
of the EL quantiles. For notational simplicity, we
will state the result only for the bivariate case.
Let $\xi_r$ be the population quantile at some level $\alpha_r$
of the $r$th population in the DRM.
We similarly define $\xi_s$ at some level $\alpha_s$.
The exact levels $\alpha_r$ and $\alpha_s$ are not
important.

%
\begin{theorem}
\label{jointdistofquan}
Assume that the conditions in Theorem~\ref{bahadur} hold for $\xi_r$
and $\xi_s$.
The centralized EL quantile under the DRM assumption
\[
\sqrt{n}(\hat\xi_{r} - \xi_{r}, \hat\xi_{s} -
\xi_{s})
\]
is asymptotically bivariate normal with mean 0 and covariance matrix
%
\begin{equation}
\label{sigEL} \Sigma_{\EL} =\pmatrix{
\omega_{rr}(\xi_r, \xi_r)/g_r^2(
\xi_r) & \omega_{rs}(\xi_r,
\xi_s)/\bigl\{g_r(\xi_r)g_s(
\xi_s)\bigr\}
\vspace*{2pt}\cr
\omega_{rs}(\xi_r, \xi_s)/\bigl
\{g_r(\xi_r)g_s(\xi_s)\bigr\}
& \omega_{ss}(\xi_s, \xi_s)/g_s^2(
\xi_s)
}.
\end{equation}
\end{theorem}

The above result does not restrict the selection of the two
populations or the levels of the quantiles.
It can be used to
conveniently obtain the limiting distributions of
smooth functions of the EL quantiles.

\section{Efficiency comparison}\label{sec4}

The EL quantiles are constructed by pooling information
from $m+1$ independent random samples.
We trust that they are more efficient than
empirical quantiles (hereafter EM) based on single samples.
A rigorous proof of this intuitive claim is not
simple.

Let $\breve{G}_{r}(x)$ be the empirical distribution function
based solely on the $r$th sample.
As processes indexed by $x$,
$
\sqrt{n}\{\breve{G}_r(x) - G_r(x)\}
$,
$r=0, 1, \ldots, m$, are independent and each
converges in distribution to a Gaussian process
with covariance function $\sigma_{rr}(x, x)$.
Let $\breve{\xi}_{r}$ be the EM quantiles of
the $r$th population at level $\alpha_r$.
Based on the classical Bahadur presentation, with any number of
choices in $r$ and $\alpha_r$, $\{ \sqrt{n}(\breve{\xi}_{r} - \xi_{r})\}$
are jointly asymptotically multivariate normal with mean 0.
In the bivariate case, the asymptotic covariance matrix of
$
\sqrt{n}(\breve{\xi}_{r} - \xi_{r}, \breve{
\xi}_{s} - \xi_{s})
$
is given by
\begin{eqnarray*}
\Sigma_{\EM} =\pmatrix{ \sigma_{rr}(
\xi_r,\xi_r) / \bigl\{ g_r^2 (
\xi_r ) \bigr\} & \sigma_{rs}(\xi_r,
\xi_s) / \bigl\{ g_r ( \xi_r )
g_s ( \xi_s ) \bigr\}
\vspace*{2pt}\cr
\sigma_{sr}(\xi_s,\xi_r) / \bigl\{
g_r ( \xi_r ) g_s ( \xi_s )
\bigr\} & \sigma_{ss}(\xi_s,\xi_s)/ \bigl\{
g_s^2(\xi_s) \bigr\}},
\end{eqnarray*}
where $\sigma_{rs}(x,y)$ was given in
Theorem~\ref{asyofdistest}.

Since the EL and EM quantiles are
asymptotically unbiased, the efficiency comparison
reduces to a comparison of two
asymptotic covariance matrices.
The following result generalizes
Corollary 4.3 of \citet{Zha00}.

%
\begin{theorem}
\label{matdiff}
For any pair of integers $0\leq r, s\leq m$
and any quantile levels $\alpha_r$ and $\alpha_s$, we have
$\Sigma_{\EM} - \Sigma_{\EL} \geq0$. This
conclusion remains true
for any number of quantiles.
\end{theorem}


\section{Inferences on functions of quantiles}\label{secquancom}

In applications such as the wood project,
we are interested in the
size of $\xi_r$, $\xi_r - \xi_s$, etc.
for various choices of $r$ and $s$ and various levels.
Two scenarios are of particular interest.
(A) For a specific wood product in a given year, is its quality
index above or below the industrial standard?
(B) How different are the quality indices for wood products
produced in two different years, mills or regions?

(A) and (B) can be addressed through
hypothesis tests or the construction of confidence
intervals. With the asymptotic normality and favorable
efficiency properties of the EL quantiles, the
task is simple. We must
find a consistent estimate of $\var(\hat\xi_r)$ and
construct approximate $100(1-\alpha)\%$ confidence intervals
as
$
\hat\xi_{r} \pm z_{\alpha/2} \sqrt{ \hat{\var}(\hat
\xi_r) },
$
where $z_{{\alpha/2}}$ denotes the $(1-{\alpha}/{2})$th
quantile of the standard normal distribution.
Similarly, approximate confidence intervals for $\xi_r - \xi_s$
are
$
(\hat\xi_{r} - \hat\xi_{s}) \pm z_{\alpha/2} \sqrt{
\hat{\var}(\hat\xi_{r} - \hat\xi_{s}) }.
$
In both cases, we need
effective and consistent estimates of $\var(\hat\xi_r)$
and $\var(\hat\xi_{r} - \hat\xi_{s})$.

With the help of (\ref{sigEL}),
plug-in consistent variance estimators can
easily be constructed.
Two necessary ingredients are consistent
estimations of $\omega_{rs}(x, y)$ and $g_r(x)$.
Although $\hat G_r(x)$ is discrete, the
idea of kernel density estimation can be used
to produce a density estimate.
Let $K(\cdot) \geq0$ be a commonly used kernel function
such that $\int K(x) \,dx = 1$, $\int x K(x) \,dx =0$ and
$\int x^2 K(x) \,dx < \infty$.
For some bandwidth $b> 0$, let
$
K_b(x) = (1/b) K( x/b).
$
Then a kernel estimate of $g_r(x)$ is given by
\[
\hat g_r(x) = \int K_b(x - y) \,d\hat
G_r(y).
\]
In the simulation study, we set $K(x)$ to the standard normal density function.
We chose the bandwidth $b$ according to the rule of thumb of
\citet{Deh77} and \citet{Sil86},
\[
b = 1.06 n^{-1/5} \min\{ \hat\sigma, \hat R/1.34\}.
\]
The above formula is designed for the situation
where the density function is estimated based on
$n$ independent and identically distributed observations.
In our simulation, we regard the fitted $\hat G_r(x)$ as a
nonrandom distribution function, and compute the
standard deviation and inter-quartile range of this distribution as
$\hat\sigma$ and $\hat R$.

The analytical form of $\omega_{rs}(x, y)$ contains many terms,
but it is straightforward to estimate them consistently and sensibly.
Let
\[
\hat a_{rs}(x) = \int_{-\infty}^x \bigl\{
\delta_{rs} h_r(t; \hat{\btheta}) - h_r(t;
\hat{\btheta}) h_s(x; \hat{\btheta}) \bigr\} h(t; \hat{\btheta}) \,d
\hat G_0(t),
\]
and we form $\hat{\bB}_r(x)$ and $\hat{\bW}$ via
\begin{eqnarray*}
\hat{\bB}_{r,s}(x) &=& \int_{-\infty}^x
\bigl\{ \delta_{rs} h_r(t; \hat{\btheta}) -
h_r(t; \hat{\btheta}) h_s(x; \hat{\btheta}) \bigr\}
\bq(t) h(t; \hat{\btheta}) \,d \hat G_0(t),
\\
\hat{\bW}_{rs} &=& \int_{-\infty}^{\infty} \bigl
\{ \delta_{rs} h_r(t; \hat{\btheta}) - h_r(t;
\hat{\btheta}) h_s(x; \hat{\btheta}) \bigr\} \bq(t)
\bq^{\tau}(t) h(t; \hat{\btheta}) \,d \hat G_0(t).
\end{eqnarray*}
We then form a consistent estimator of $\omega_{rs}(\xi_r, \xi_s)$ as
\[
\hat\omega_{rs}(\xi_r, \xi_s) =
\rho_r^{-1} \delta_{rs} ( \alpha_r -
\alpha_r \alpha_s) - (\rho_r
\rho_s)^{-1} \bigl\{ \hat a_{rs}( \hat
\xi_r \wedge\hat\xi_s) - \hat{\bB}_r^\tau(
\hat\xi_r) \hat{\bW}^{-1} \hat{\bB}_s(\hat
\xi_s) \bigr\},
\]
where we have used the facts that
$\delta_{rs}G_r(\xi_r \wedge\xi_s) = \delta_{rs} \alpha_r $
and $G_r(\xi_r) = \alpha_r$.


\section{Simulation study}\label{sec6}

We now examine the finite-sample performance of the
inference procedures via simulation.
Are the EL quantiles $\hat\xi$ more
efficient than the EM quantiles $\breve\xi$?
The simulation studies shed light on how large the sample must be
before the asymptotic result applies.
We analyze
data sets generated from several sets of populations, which
are divided into
two groups: those that satisfy the DRM assumption and those that do not.

\subsection{Populations satisfying DRM assumptions}\label{sec6.1}
Recall that the Gamma and normal distribution families are
special DRMs. We choose two sets of distributions
from these families with the parameter values
specified in Table~\ref{table1}.
For the Gamma distributions, the first parameter is
the degrees of freedom and the second is the scale.
Therefore, the expectation of the first population is $6 \times1.5$.
The parameters for the normal distribution are the mean and
variance.
The populations have
similar means and variances to those seen in applications.

%
\begin{table}
\caption{Parameters under DRM}\label{table1}
\begin{tabular*}{\textwidth}{@{\extracolsep{\fill}}lcccccc@{}}
\hline
\multicolumn{6}{c}{\textbf{Distributions}} & \multicolumn
{1}{c@{}}{$\bolds{\bq^{\tau}(x)}$} \\
\hline
$\Gamma$(6, 1.5) & $\Gamma$(6, 1.4) & $\Gamma$(7, 1.3)
& $\Gamma$(7, 1.2) & $\Gamma$(8, 1.1) & $\Gamma$(8, 1.0)
& ($1, x, \log(x)$) \\
$N(18, 4)$ & $N(18, 9)$ & $N(20, 6)$
& $N(20, 9)$ & $N(22, 8)$ & $N(22, 10)$ &($1, x, x^2$) \\
\hline
\end{tabular*}
\end{table}

The simulations were carried out with
$n_r=50$ and 2000 repetitions.
We examine the performance of $\hat G_r(x)$
and $\breve G_r(x)$ for $x$ set to the quantile
levels $\alpha= 5\%$, 10\%, 50\%, 90\% and 95\%.
We computed the relative bias,
the asymptotic variance and the simulated variance of the
EL estimator $\hat G_r(x)$.
The EM estimator $\breve G_r(x)$ is unbiased, and
its asymptotic variance is $\alpha(1-\alpha)/\rho_r$.
For ease of comparison,
we report the ratios of the EM and EL asymptotic variances
and the ratios of their simulated variances.
We also report the ratios of the
mean estimated variances of $\hat G_r(x)$
and their corresponding asymptotic variances.
The results are presented in Table~\ref{cdf_comp_gamma_50}.

%
\begin{table}
\def\arraystretch{0.9}
\tabcolsep=0pt
\caption{EL and EM distribution estimates.
$\hat G$: EL estimate;
$\breve G$: EM estimate;
$\sigma^2$: asymptotic variance;
$\hat\sigma^2$: average of variance estimate;
${\mathrm B}$: bias of $\sqrt{\alpha(1-\alpha)}$ as percentage;
$n_k =50$}\label{cdf_comp_gamma_50}
\begin{tabular*}{\textwidth}{@{\extracolsep{\fill}}lccccccc@{}}
\hline
& \multicolumn{1}{c}{$\bolds{\alpha}$} & \multicolumn{1}{c}{$\bolds
{\sigma^2(\hat G)}$} & \multicolumn{1}{c}{$ \bolds{\sigma^2(\breve
G)/\sigma^2(\hat G)} $} &
\multicolumn{1}{c}{$\bolds{\var( \breve G)/\var( \hat G)}$} &
\multicolumn{1}{c}{$\bolds{\hat\sigma^2(\hat G)/\var(\hat G)}$} &
\multicolumn{1}{c}{$\bolds{ {\mathrm B}(\hat G)}$} &
\multicolumn{1}{c@{}}{$ \bolds{{\mathrm B}(\breve G)} $} \\
\hline
$\Gamma(6, 1.5) $ & 0.05& 0.176 & 1.62& 1.48&1.03 & $-$0.46& $-$4.02 \\
& 0.10& 0.378 & 1.43& 1.38&1.00 & $-$0.28& $-$3.24 \\
& 0.50& 1.055 & 1.42& 1.33&0.98 & $-$0.65& $-$1.92 \\
& 0.90& 0.382 & 1.41& 1.38&1.01 & $-$0.83& $-$3.68 \\
& 0.95& 0.178 & 1.60& 1.58&1.12 & $-$0.45& $-$4.67 \\[3pt]
$\Gamma(6, 1.4)$& 0.05& 0.176 & 1.62& 1.49&1.08 & $-$1.15& $-$4.37 \\
& 0.10& 0.374 & 1.44& 1.39&1.03 & $-$0.71& $-$3.33 \\
& 0.50& 1.031 & 1.45& 1.37&1.00 & $-$0.31& $-$2.03 \\
& 0.90& 0.370 & 1.46& 1.39&1.04 & $-$0.68& $-$3.20 \\
& 0.95& 0.172 & 1.66& 1.66&1.13 & $-$0.74& $-$4.39 \\[3pt]
$\Gamma(7, 1.3)$& 0.05& 0.170 & 1.67& 1.50&1.00 & $-$0.27& $-$4.18 \\
& 0.10& 0.368 & 1.47& 1.49&1.04 & $-$0.20& $-$2.74 \\
& 0.50& 1.034 & 1.45& 1.44&0.98 & $-$0.29& $-$1.74 \\
& 0.90& 0.394 & 1.37& 1.39&1.01 & $-$0.24& $-$3.12 \\
& 0.95& 0.186 & 1.53& 1.60&1.07 & $-$0.29& $-$4.33 \\[3pt]
$\Gamma(7, 1.2)$& 0.05& 0.171 & 1.67& 1.56&1.06 & $-$0.86& $-$4.14 \\
& 0.10& 0.369 & 1.46& 1.39&1.05 & $-$0.78& $-$3.38 \\
& 0.50& 1.065 & 1.41& 1.32&0.97 & $-$0.51& $-$2.33 \\
& 0.90& 0.371 & 1.46& 1.42&0.99 & $-$0.81& $-$3.57 \\
& 0.95& 0.169 & 1.68& 1.60&0.99 & $-$1.13& $-$4.68 \\[3pt]
$\Gamma(8, 1.1)$ & 0.05& 0.172 & 1.65& 1.54&1.06 & $-$0.72& $-$4.26 \\
& 0.10& 0.372 & 1.45& 1.44&1.05 & $-$0.67& $-$3.27 \\
& 0.50& 1.055 & 1.42& 1.43&1.01 & $-$0.25& $-$1.53 \\
& 0.90& 0.373 & 1.45& 1.47&1.01 & $-$0.19& $-$3.56 \\
& 0.95& 0.176 & 1.62& 1.68&1.08 & $-$0.37& $-$4.15 \\[3pt]
$\Gamma(8, 1.0)$ & 0.05& 0.180 & 1.59& 1.55&1.07 & $-$1.04& $-$4.22 \\
& 0.10& 0.379 & 1.42& 1.41&0.98 & $-$0.81& $-$3.44 \\
& 0.50& 1.041 & 1.44& 1.33&0.95 & \phantom{$-$}0.28 & $-$1.63 \\
& 0.90& 0.364 & 1.48& 1.52&1.09 & \phantom{$-$}0.01 & $-$2.99 \\
& 0.95& 0.164 & 1.73& 1.64&1.09 & $-$0.64& $-$4.39 \\
\hline
\end{tabular*}
\end{table}

%
\begin{table}
\def\arraystretch{0.9}
\caption{EL and EM quantiles.
$\xi$: true quantile;
$\hat\xi$: EL quantile;
$\breve\xi$: EM quantile;
$\sigma^2$: asymptotic variance;
$\hat\sigma^2$: average of variance estimate;
${\mathrm B}$: relative bias as percentage;
$n_k =50$}\label{quan_comp_gamma_50}
\begin{tabular*}{\textwidth}{@{\extracolsep{\fill
}}lcd{2.2}d{3.2}d{3.2}d{3.2}cd{2.2}d{2.2}@{}}
\hline
& \multicolumn{1}{c}{$\bolds{\alpha}$} & \multicolumn{1}{c}{$\bolds{\xi
}$}& \multicolumn{1}{c}{$\bolds{\sigma^2(\hat\xi)}$} &
\multicolumn{1}{c}{$\bolds{ \var( \hat\xi)} $} &
\multicolumn{1}{c}{$\bolds{\hat{\sigma}^2( \hat\xi) }$} &
\multicolumn{1}{c}{$\bolds{\var( \breve\xi)/\var( \hat\xi)}$} &
\multicolumn{1}{c}{$\bolds{ \mathrm{ B}(\hat\xi) }$} &
\multicolumn{1}{c@{}}{$\bolds{\mathrm{ B}(\breve\xi)}$} \\
\hline
$\Gamma(6, 1.5)$ & 0.05& 3.92& 71.31& 70.38& 67.50& 1.51& 2.19 & 1.59\\
& 0.10& 4.73& 69.23& 68.28& 73.79& 1.34& 1.08 & 0.97\\
& 0.50& 8.51& 83.70& 88.50& 94.46& 1.27& 0.45 & 0.15\\
& 0.90& 13.91& 298.78& 292.15& 287.19& 1.35& 0.23 & 0.19\\
& 0.95& 15.77& 473.53& 468.31& 469.48& 1.54&-0.25 &-0.32\\[3pt]
$\Gamma(6, 1.4)$
& 0.05& 3.66& 62.34& 62.42& 59.15& 1.55& 2.84 & 1.94\\
& 0.10& 4.41& 59.63& 60.02& 64.02& 1.31& 1.55 & 0.88\\
& 0.50& 7.94& 71.30& 75.64& 81.58& 1.39& 0.27 & 0.28\\
& 0.90& 12.98& 252.32& 248.04& 248.28& 1.37& 0.11 &-0.05\\
& 0.95& 14.72& 397.73& 387.44& 399.88& 1.51&-0.05 &-0.41\\[3pt]
$\Gamma(7, 1.3)$
& 0.05& 4.27& 67.35& 66.55& 62.81& 1.59& 1.83 & 1.07\\
& 0.10& 5.06& 63.93& 64.70& 68.40& 1.33& 0.88 & 0.50\\
& 0.50& 8.67& 72.61& 76.18& 82.87& 1.37& 0.30 & 0.04\\
& 0.90& 13.69& 260.55& 254.74& 241.91& 1.30&-0.08 &-0.12\\
& 0.95& 15.40& 413.83& 414.73& 395.77& 1.52&-0.22 &-0.11 \\[3pt]
$\Gamma(7, 1.2)$
& 0.05& 3.94& 57.63& 59.87& 54.10& 1.54& 2.32 & 1.11\\
& 0.10& 4.67& 54.66& 56.15& 58.62& 1.36& 1.32 & 0.75\\
& 0.50& 8.00& 63.76& 65.99& 70.68& 1.32& 0.30 & 0.23\\
& 0.90& 12.64& 208.95& 229.25& 206.39& 1.32& 0.24 & 0.12\\
& 0.95& 14.21& 320.80& 340.92& 319.98& 1.61& 0.13 &-0.13\\[3pt]
$\Gamma(8, 1.1)$
& 0.05& 4.38& 60.62& 56.57& 54.22& 1.66& 1.95 & 1.07\\
& 0.10& 5.12& 56.23& 54.25& 58.29& 1.33& 1.12 & 0.63\\
& 0.50& 8.44& 61.08& 60.83& 67.63& 1.40& 0.23 &-0.10\\
& 0.90& 12.95& 195.71& 198.99& 188.42& 1.36&-0.15 &-0.03\\
& 0.95& 14.46& 309.32& 323.05& 295.76& 1.47&-0.24 &-0.55\\[3pt]
$\Gamma(8, 1.0)$& 0.05& 3.98& 52.12& 50.42& 44.41& 1.55& 2.48 & 1.22\\
& 0.10& 4.66& 47.38& 47.27& 47.85& 1.31& 1.35 & 0.88\\
& 0.50& 7.67& 49.81& 53.38& 55.67& 1.35&-0.02 & 0.00\\
& 0.90& 11.77& 157.93& 162.29& 155.97& 1.33&-0.22 &-0.30\\
& 0.95& 13.15& 238.30& 247.86& 242.16& 1.55&-0.14 &-0.32\\
\hline
\end{tabular*}
\end{table}

There is an efficiency gain in the
range of 40\% to 70\% for the EL estimators in terms
of both the theoretical and simulated variances.
The variances of the EL estimators are estimated
accurately: in the $\hat\sigma^2(\hat G)/\sigma^2(\hat G)$
column all the entries are close to 1.
Finally, the relative biases $B(\hat G)$
and $B(\breve G)$ are both small.

We now turn to investigating the performance of the EL and EM quantiles
for both point and interval estimations.
The quantile of a discrete distribution $G(x)$ is not a smooth function,
and this puts the EM quantile at a disadvantage.
To ensure that the EL quantile had a strong competitor,
we modified the EM quantile. We replaced $\breve G(x)$ by
$\breve G_r(x) - (2n_r)^{-1}$, and we used linear interpolation
to calculate this quantile.
These modifications do not alter the first-order asymptotics.
We continue to use the notation $\hat\xi_{r}$ and $\breve\xi_{r}$
for the EL and EM quantiles after these modifications.

The simulation results for the quantile estimates are
given in Table~\ref{quan_comp_gamma_50}.
The simulated EL variances $\var(\hat\xi)$
and the mean estimated EL variances $\hat\sigma^2(\hat\xi)$
are both close to the asymptotic variances
$\sigma^2(\hat\xi)$.
The results support the asymptotic theory and
the viability of the EL variance estimator.
The ratio $\var(\breve\xi)/\var(\hat\xi)$ is based on simulated variances
and ranges between 1.20 and 1.60.
These results indicate an efficiency gain of between 20\% and 60\%
in the EL quantiles.
Finally, the relative biases $B(\hat\xi)$ and $B(\breve\xi)$ are low and
within $\pm3$\%.

The simulation results for the interval estimates of the
quantiles and quantile
differences are given in Table~\ref{length-cov-gam}.
We compute the average lengths and coverage probabilities of
the EL and EM confidence intervals at the $95\%$ level.
The coverage probabilities of the EL intervals
are almost always closer to the nominal $95\%$.
This advantage is more obvious for the upper-tail quantiles
(such as the 95\% quantile).
Often, the coverage gains of the EL intervals reach 5\%,
and these intervals are 10\% to 20\% shorter.

%
\begin{table}
\def\arraystretch{0.9}
\tabcolsep=0pt
\caption{Confidence intervals for quantile and quantile differences.
Nominal level: $95\%$;
$\Gamma(6, 1.5) -\Gamma(6, 1.4)$: differences of
$\Gamma(6, 1.5)$ and $\Gamma(6, 1.4)$ quantiles\break
at the given level; $n_k=50$}\label{length-cov-gam}
{\fontsize{8.5}{10.5}{\selectfont  \begin{tabular*}{\textwidth}{@{\extracolsep{\fill
}}lcd{2.2}d{2.2}d{2.2}d{2.2}d{2.2}d{2.2}d{2.2}d{2.2}d{2.2}d{2.2}@{}}
\hline
& & \multicolumn{5}{c}{\textbf{EL}} & \multicolumn{5}{c@{}}{\textbf{EM}}
\\[-6pt]
& & \multicolumn{5}{c}{\hrulefill} & \multicolumn{5}{c@{}}{\hrulefill}
\\
& \multicolumn{1}{c}{$\bolds{\alpha}$} & \multicolumn{1}{c}{$\bolds{5\%
}$}& \multicolumn{1}{c}{$\bolds{10\%}$}&
\multicolumn{1}{c}{$\bolds{50\%}$} & \multicolumn{1}{c}{$\bolds{90\%
}$}& \multicolumn{1}{c}{$\bolds{95\%}$} &
\multicolumn{1}{c}{$\bolds{5\%}$} &\multicolumn{1}{c}{$\bolds{10\%}$}&
\multicolumn{1}{c}{$\bolds{50\%}$}&
\multicolumn{1}{c}{$\bolds{90\%}$} & \multicolumn{1}{c}{$\bolds{95\%}$}
\\
\hline
$\Gamma(6, 1.5)$ &length & 1.83& 1.93& 2.18& 3.77& 4.73& 2.33& 2.42&
2.66& 4.23& 5.20 \\
&coverage& 93.0& 94.8& 95.2& 90.7& 89.5& 92.3& 96.5& 95.9& 89.1& 84.9 \\[3pt]
$\Gamma(6, 1.4)$ &length & 1.71& 1.79& 2.03& 3.49& 4.36& 2.16& 2.24&
2.47& 3.98& 4.87 \\
&coverage& 91.8& 94.3& 95.0& 90.9& 88.7& 91.7& 95.0& 96.3& 89.0& 84.3 \\[3pt]
$\Gamma(7, 1.3)$ &length & 1.76& 1.85& 2.03& 3.41& 4.27& 2.25& 2.31&
2.48& 3.88& 4.78 \\
&coverage& 92.0& 93.7& 95.5& 92.2& 90.1& 90.0& 95.0& 96.8& 90.1& 85.1 \\[3pt]
$\Gamma(7, 1.2)$ &length & 1.63& 1.71& 1.88& 3.16& 3.92& 2.07& 2.13&
2.31& 3.62& 4.45 \\
&coverage& 91.9& 94.3& 95.3& 92.5& 91.3& 91.3& 95.5& 97.4& 89.4& 85.1 \\[3pt]
$\Gamma(8, 1.1)$ &length & 1.64& 1.71& 1.84& 3.02& 3.73& 2.09& 2.12&
2.26& 3.50& 4.24 \\
&coverage& 91.5& 93.9& 95.4& 92.0& 90.1& 91.5& 95.0& 95.9& 90.6& 85.7 \\[3pt]
$\Gamma(8, 1.0)$ &length & 1.50& 1.56& 1.67& 2.74& 3.36& 1.91& 1.93&
2.05& 3.15& 3.82 \\
&coverage& 90.1& 93.4& 95.3& 92.5& 91.7& 89.3& 94.9& 96.6& 90.0& 85.1
\\[6pt]
$\Gamma(6, 1.5) -\Gamma(6, 1.4)$ &length & 2.41& 2.58& 2.92& 5.08&
6.33& 3.23& 3.32& 3.65& 5.95& 7.38 \\
&coverage& 95.0& 95.5& 95.5& 93.8& 94.5& 94.5& 97.3& 96.0& 91.8& 89.5 \\[3pt]
$\Gamma(6, 1.5) -\Gamma(7, 1.3)$ &length & 2.47& 2.62& 2.91& 5.00&
6.22& 3.29& 3.37& 3.66& 5.87& 7.33 \\
&coverage& 95.1& 96.1& 95.3& 94.8& 95.1& 94.1& 96.8& 96.2& 93.2& 89.9 \\[3pt]
$\Gamma(6, 1.5) -\Gamma(7, 1.2)$ &length & 2.37& 2.53& 2.82& 4.89&
6.09& 3.17& 3.25& 3.54& 5.69& 7.09 \\
&coverage& 93.9& 95.5& 95.8& 94.1& 93.7& 94.2& 96.2& 96.5& 93.4& 90.2 \\[3pt]
$\Gamma(6, 1.5) -\Gamma(8, 1.1)$ &length & 2.40& 2.54& 2.79& 4.78&
5.95& 3.18& 3.25& 3.51& 5.61& 6.95 \\
&coverage& 94.2& 95.8& 96.3& 94.7& 93.7& 94.4& 96.6& 96.1& 93.2& 90.0 \\[3pt]
$\Gamma(6, 1.5) -\Gamma(8, 1.0)$ &length & 2.29& 2.43& 2.71& 4.68&
5.84& 3.06& 3.12& 3.38& 5.39& 6.67 \\
&coverage& 93.5& 95.2& 95.9& 92.3& 91.9& 94.1& 96.4& 96.7& 92.7& 89.0 \\
\hline
\end{tabular*}   }}\vspace*{6pt}
\end{table}

We also conducted simulations for the second group of
populations, as shown in Table~\ref{table1}, and for
$n_r = 100$. The results are similar and omitted.
In conclusion, the EL approach is superior when the model
assumptions are satisfied.

\subsection{Performance when model is misspecified}\label{sec6.2}
What happens to the EL approach
when the model is misspecified?
\citet{FokKai06} quantified the effect of choosing an
incorrect linear form of $\bq(x)$.
In general, both the point estimation and the hypothesis tests
are adversely affected when the model is misspecified.
These findings may have motivated the model selection approach
in \citet{Fok07}. That is, instead of pre-specifying a known
$\bq(x)$, one may select $\bq(x)$ as a linear combination of
a rich class of functions. For instance, let $\bq(x) = \{1, \log x,
x^{0.5}, x, x^{1.5}, x^2\}^\tau$.
The most appropriate $\bq(x)$ is then determined by
selecting a subvector of the current $\bq(x)$.
Hence, the classical model selection approaches can be used.

Following this lead, we provide a limited study of
the impact of misspecification on the quantile estimations.
For this purpose, we simulated random samples
from a number of Gamma distributions, Weibull distributions,
denoted $W(\cdot, \cdot)$, and normal distributions, as shown in Table~\ref{table2}.
These populations are chosen to have similar means and variances.
We obtained the EL quantile estimates as if they satisfy DRM for some
pre-specified but wrong $\bq(x)$.

%
\begin{table}
\caption{Parameters for non-DRM}
\label{table2}
\begin{tabular*}{\textwidth}{@{\extracolsep{\fill}}lccccc@{}}
\hline
$\Gamma$(16, 0.6) & $\Gamma$(19, 0.5)
& $N(9, 5)$ & $N(9.6, 5.6)$ & $W(10, 4.5)$ & $W(11, 5)$ \\
$\Gamma$(16, 0.6) & $\Gamma$(19, 0.5) & $\Gamma$(17.5, 0.5)
& $W(10.5, 4.5)$ & $W(10, 4.5)$ & $W(11, 5)$ \\
\hline
\end{tabular*}
\end{table}

%
\begin{table}[b]
\def\arraystretch{0.9}
\caption{EL and EM quantiles under model mis-specification.
$\xi$: true quantile;
$\hat\xi$: EL quantile;
$\breve\xi$: EM quantile;
$\hat\sigma^2$: average of variance estimate;
${\mathrm B}$: relative bias as percentage;
$n_k =50$}\label{quan_comp_misn_50}
\begin{tabular*}{\textwidth}{@{\extracolsep{\fill
}}lcd{2.3}d{2.3}ccd{2.2}d{2.2}@{}}
\hline
& \multicolumn{1}{c}{$\bolds{\alpha}$} & \multicolumn{1}{c}{$\bolds{\xi
} $}& \multicolumn{1}{c}{$\bolds{ \mse( \hat\xi)} $} &
\multicolumn{1}{c}{$\bolds{ \mse( \breve\xi)/\mse( \hat\xi)}$} &
\multicolumn{1}{c}{$\bolds{\hat\sigma(\hat\xi)/\mse( \hat\xi)}$} &
\multicolumn{1}{c}{$\bolds{ {\mathrm B}(\hat\xi) }$} &
\multicolumn{1}{c@{}}{$ \bolds{ {\mathrm B}(\breve\xi) }$} \\
\hline
$\Gamma(16, 0.6)$ & 0.05 & 6.022& 61.06& 1.23& 0.83& 1.20& 0.56 \\
& 0.10 & 6.681& 49.19& 1.19& 0.95& 0.48& 0.07 \\
& 0.50 & 9.401& 44.14& 1.16& 1.13& -0.16& -0.18 \\
& 0.90 & 12.775& 125.38& 1.10& 0.94& 0.24& -0.04 \\
& 0.95 & 13.858& 210.07& 1.17& 0.92& 0.05& -0.27 \\[2ex]
$\Gamma(19, 0.5)$
& 0.05 & 6.221& 57.81& 1.22& 0.78& 1.29& 0.71 \\
& 0.10 & 6.836& 44.47& 1.14& 0.90& 0.77& 0.33 \\
& 0.50 & 9.334& 34.15& 1.20& 1.17& 0.07& 0.10 \\
& 0.90 & 12.378& 97.34& 1.14& 0.97& 0.17& -0.00 \\
& 0.95 & 13.346& 151.68& 1.20& 0.98& 0.06& -0.40 \\[2ex]
$N(9, 5)$
& 0.05 & 5.322& 103.06& 1.23& 0.85& 1.43& 0.77 \\
& 0.10 & 6.134& 73.47& 1.08& 0.94& 0.17& 0.23 \\
& 0.50 & 9.000& 37.86& 1.20& 1.14& 0.21& -0.03 \\
& 0.90 & 11.866& 68.81& 1.18& 0.92& -0.26& -0.17 \\
& 0.95 & 12.678& 97.39& 1.30& 0.89& -0.39& -0.39 \\[2ex]
$N(9.6, 5.6)$
& 0.05 & 5.708& 116.53& 1.25& 0.89& 0.79& 0.49 \\
& 0.10 & 6.567& 87.51& 1.14& 0.90& 0.18& 0.13 \\
& 0.50 & 9.600& 43.71& 1.20& 1.11& 0.37& 0.27 \\
& 0.90 & 12.633& 76.92& 1.22& 0.96& -0.06& -0.13 \\
& 0.95 & 13.492& 112.45& 1.25& 0.91& -0.26& -0.30 \\[2ex]
$W(10, 4.5)$
& 0.05 & 3.344& 12.09& 1.01& 0.85& 0.89& 0.48 \\
& 0.10 & 3.593& 6.62& 1.01& 0.97& 0.58& 0.20 \\
& 0.50 & 4.338& 2.08& 1.10& 1.11& -0.09& -0.01 \\
& 0.90 & 4.891& 1.91& 1.17& 1.14& -0.01& -0.11 \\
& 0.95 & 5.022& 2.17& 1.39& 1.17& 0.13& -0.18 \\[2ex]
$W(11, 5.0)$
& 0.05 & 3.817& 11.89& 1.15& 0.94& 0.79& 0.25 \\
& 0.10 & 4.075& 6.67& 1.11& 0.94& 0.56& 0.01 \\
& 0.50 & 4.836& 1.99& 1.16& 1.14& -0.09& -0.03 \\
& 0.90 & 5.394& 2.15& 1.03& 1.18& -0.11& -0.07 \\
& 0.95 & 5.524& 2.63& 1.14& 1.18& 0.05& -0.08 \\
\hline
\end{tabular*}
\end{table}

As a trade-off between model interpretation and parsimony,
we choose $\bq(x) = (1, x, \log(1+|x|), \sqrt{|x|})^{\tau}$.
The remaining settings are the same as before.
In Table~\ref{quan_comp_misn_50}, we report only
the biases and mean square errors (mse) of the EL and EM
quantiles.
The EL quantiles are still uniformly more efficient with
the efficiency gains usually above 15\%.
The variance estimators remain accurate, and
the relative biases $B(\hat\xi)$ and $B(\breve\xi)$ are
still negligible.

The simulation results for the interval estimates of the
quantiles and quantile differences are given in Table~\ref{length-cov-mis}.
We compute the average lengths and the coverage probabilities of
the EL and EM confidence intervals at the $95\%$ level.
The EL confidence intervals are not clearly better.
These intervals have better coverage probabilities for the upper quantiles
but similar or slightly inferior probabilities for the lower quantiles.
The EL intervals are always shorter,
and they are more than 10\% shorter in most cases.
The simulation results for the second set of populations are similar;
they are omitted to save space.

%
\begin{table}
\tabcolsep=0pt
\caption{Confidence intervals for quantile and quantile differences
under model mis-specification.
Nominal level $95\%$; $n_k = 50$}\label{length-cov-mis}
{\fontsize{8.5}{10.5}{\selectfont
\begin{tabular*}{\textwidth}{@{\extracolsep{\fill
}}lcd{2.2}d{2.2}d{2.2}d{2.2}d{2.2}d{2.2}d{2.2}d{2.2}d{2.2}d{2.2}@{}}
\hline
& & \multicolumn{5}{c}{\textbf{EL}} & \multicolumn{5}{c@{}}{\textbf{EM}}
\\[-6pt]
& & \multicolumn{5}{c}{\hrulefill} & \multicolumn{5}{c@{}}{\hrulefill}
\\
& \multicolumn{1}{c}{$\bolds{\alpha}$} & \multicolumn{1}{c}{$\bolds{5\%
}$}& \multicolumn{1}{c}{$\bolds{10\%}$}&
\multicolumn{1}{c}{$\bolds{50\%}$} & \multicolumn{1}{c}{$\bolds{90\%
}$}& \multicolumn{1}{c}{$\bolds{95\%}$}
&\multicolumn{1}{c}{$\bolds{5\%}$} &\multicolumn{1}{c}{$\bolds{10\%}$}&
\multicolumn{1}{c}{$\bolds{50\%}$}&
\multicolumn{1}{c}{$\bolds{90\%}$} & \multicolumn{1}{c@{}}{$\bolds{95\%
}$} \\
\hline
$\Gamma(16, 0.6)$ &length & 1.57& 1.52& 1.57& 2.36& 2.91& 1.83& 1.79&
1.77& 2.50& 3.00 \\
&coverage& 87.5& 91.8& 95.7& 91.2& 88.9& 89.0& 94.1& 96.5& 90.6& 86.4 \\[3pt]
$\Gamma(19, 0.5)$ &length & 1.48& 1.41& 1.40& 2.10& 2.63& 1.71& 1.64&
1.61& 2.25& 2.66 \\
&coverage& 87.9& 90.8& 95.5& 91.9& 89.7& 88.5& 93.0& 96.4& 91.0& 85.2 \\[3pt]
$N(9, 5)$ &length & 2.04& 1.82& 1.47& 1.77& 2.05& 2.22& 2.00& 1.67&
2.02& 2.28 \\
&coverage& 87.4& 92.0& 94.6& 91.1& 89.6& 86.2& 92.3& 96.4& 91.7& 87.6 \\[3pt]
$N(9.6, 5.6)$ &length & 2.13& 1.90& 1.54& 1.91& 2.22& 2.35& 2.12& 1.76&
2.12& 2.35 \\
&coverage& 83.5& 90.6& 96.1& 92.4& 89.3& 84.4& 91.3& 96.8& 92.3& 86.6 \\[3pt]
$W(10, 4.5)$ &length & 0.66& 0.54& 0.34& 0.33& 0.35& 0.69& 0.57& 0.37&
0.36& 0.37 \\
&coverage& 82.3& 88.8& 95.1& 96.1& 95.3& 84.7& 90.8& 96.0& 95.0& 91.6 \\[3pt]
$W(11, 5.0)$ &length & 0.68& 0.53& 0.33& 0.35& 0.38& 0.70& 0.59& 0.37&
0.36& 0.37 \\
&coverage& 83.3& 87.9& 93.8& 93.7& 92.4& 83.5& 90.9& 94.8& 94.6& 89.9 \\
[6pt]
$\Gamma(16, 0.6) - \Gamma(19, 0.5) $ & length & 2.14& 2.04& 2.08& 3.15&
3.95& 2.55& 2.45& 2.41& 3.42& 4.14 \\
& coverage& 92.8& 94.6& 95.4& 94.5& 93.3& 93.5& 95.8& 96.1& 94.2& 91.3
\\[3pt]
$\Gamma(16, 0.6) -N(9, 5) $ & length & 2.59& 2.36& 2.13& 2.96& 3.60&
2.94& 2.71& 2.45& 3.26& 3.87 \\
& coverage& 92.4& 94.3& 95.7& 93.4& 91.9& 91.4& 95.3& 96.9& 94.0& 90.9
\\[3pt]
$\Gamma(16, 0.6) -N(9.6, 5.6) $ & length & 2.65& 2.41& 2.17& 3.01&
3.66& 3.05& 2.81& 2.51& 3.33& 3.92 \\
& coverage& 90.5& 93.4& 96.5& 94.5& 93.0& 91.1& 94.7& 96.1& 94.4& 91.8
\\[3pt]
$\Gamma(16, 0.6) -W(10, 4.5) $ & length & 1.73& 1.63& 1.61& 2.39& 2.93&
1.99& 1.89& 1.81& 2.53& 3.03 \\
& coverage& 90.1& 93.3& 95.3& 91.6& 88.9& 91.8& 94.7& 96.1& 91.1& 86.1
\\[3pt]
$\Gamma(16, 0.6) -W(11, 5.0)$ & length & 1.74& 1.62& 1.60& 2.39& 2.94&
1.99& 1.89& 1.81& 2.53& 3.03 \\
& coverage& 89.3& 92.9& 95.4& 91.4& 89.0& 90.8& 94.2& 96.3& 90.8& 86.7
\\
\hline
\end{tabular*}}}
\end{table}

In conclusion, while the model misspecification has a serious impact on
the estimation of $\btheta$, as shown by \citet{FokKai06},
the quantile estimations are not as badly affected.

\section{Real-data analysis}\label{sec7}
In this section, we apply our method to lumber data.
The data come from tests conducted at an FPInnovations laboratory.
They contain the MOE and MOR measurements for lumber
produced in 2007 and in 2010
with sample sizes 98 and 282, respectively.
We analyze the MOE and MOR characteristics separately.
We regard the measurements of each index
as two independent random samples from two populations
satisfying the DRM assumption.

We use the EL approach to obtain point estimates and
confidence intervals for the quantiles
and the quantile differences between 2007 and 2010 of
each quality index.
We set $\bq(x) = (1, x, \log(1+|x|), \sqrt{|x|})^{\tau}$
as in the second simulation study.
Different choices of $\bq(x)$ do not markedly change
the quantile estimates and confidence intervals,
although they may give very different estimates of
$\btheta$.

\begin{figure}

\includegraphics{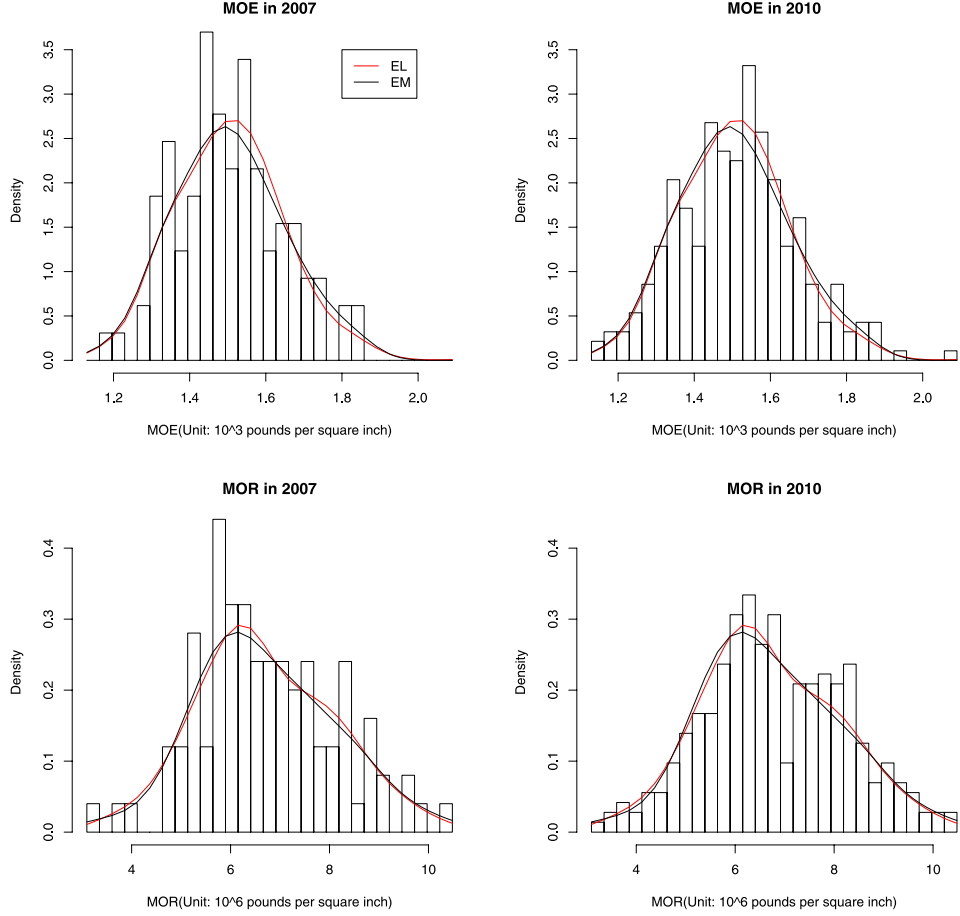}

\caption{Histograms of MOE and MOR
with EL and EM density estimates.}\label{hist-rd}
\end{figure}

\begin{figure}

\includegraphics{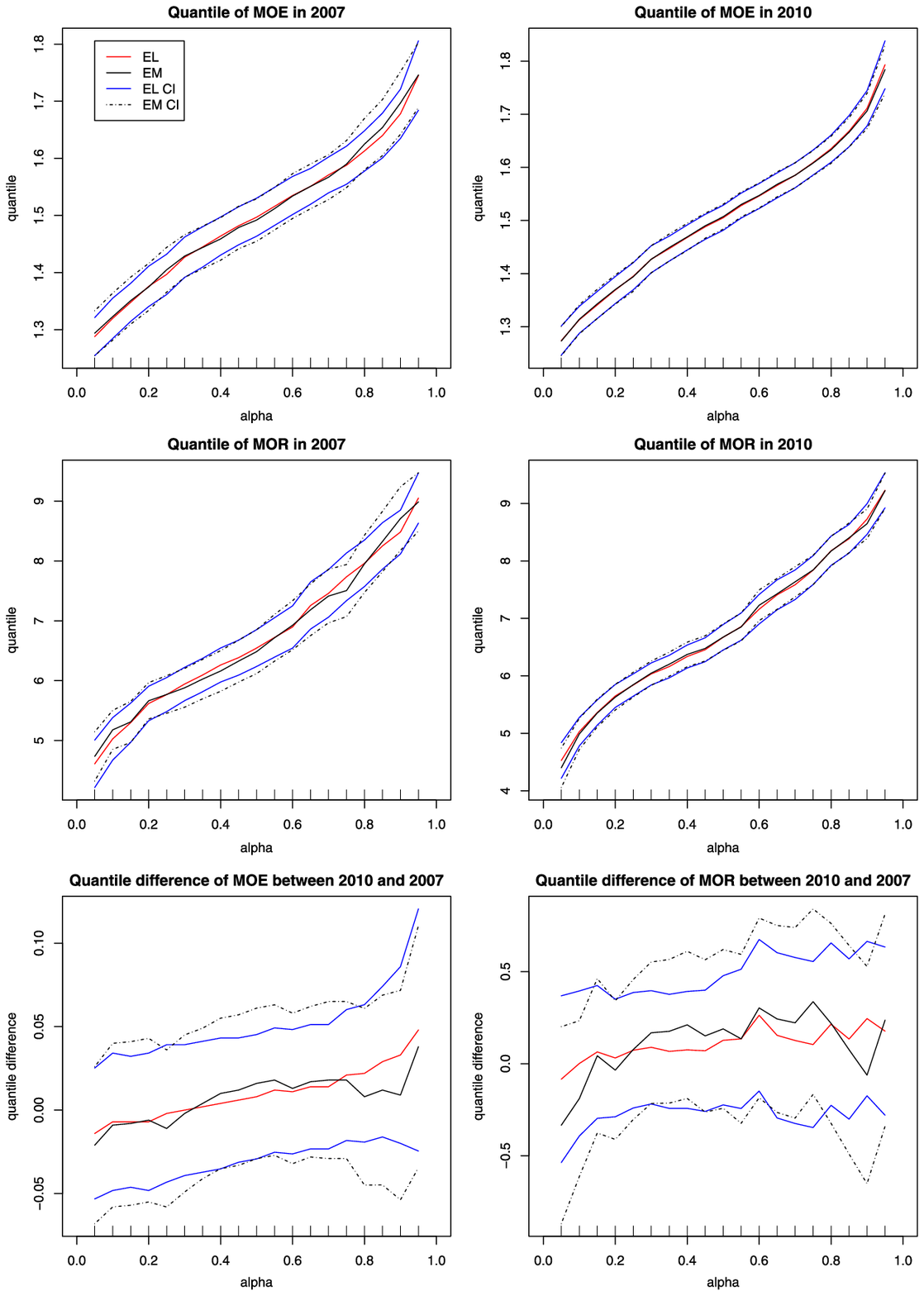}

\caption{Point estimates and confidence intervals
for quantiles and quantile differences for lumber data.}\label{ci-rd}
\end{figure}

Figure~\ref{hist-rd} presents
histograms of the MOE and MOR measurements
with the EL and EM density estimates.
We computed the EL and EM quantiles, the quantile differences,
and their 95\% confidence intervals for the
5\% level to the 95\% level in 5\% increments.
These point estimates and the
confidence limits are connected to obtain the six plots
shown in Figure~\ref{ci-rd}.
The EL quantiles and confidence limits
are much smoother than those of EM.
This phenomenon can be explained by the fact
that the EL method is designed to use information from
all the samples, which leads to less variation.
These plots do not indicate
that the EL method has sharper confidence limits.
In fact, the EL intervals are 10\% shorter than
the EM intervals for both the quantiles and the quantile differences.
In view of the simulation support for the validity of both the EL and EM
approaches, the 10\% gain is likely real, and it implies significant
cost savings in applications.
To save space, we do not include tables of the results.

%

%

\begin{appendix}
\section*{Appendix}\label{app}
\subsection{\texorpdfstring{Sketched proof of Theorem \protect\ref{thetaest}}{Sketched proof of Theorem 2.1}}\label{appa1}
This theorem mostly summarizes and extends the results in
\citet{QinZha97}, \citet{Foketal01} and \citet{Zha02}.
To enable readers to understand the other proofs, we provide
some necessary details.
Interested readers can contact the authors for more detailed derivations.

%
\begin{lemma}
\label{taylorexpan} Assume the conditions of Theorem~\ref{thetaest}.
For any $\btheta$ such that $\btheta= \btheta^* + o(n^{-1/3})$, we have
\[
\ell_n (\btheta) - \ell_n \bigl(\btheta^* \bigr) =
\bigl(\btheta- \btheta^* \bigr)^{\tau} \bZ_n -
\frac{n}{2} \bigl(\btheta- \btheta^* \bigr)^{\tau} \bW\bigl(\btheta
- \btheta^* \bigr) + o_p(1),
\]
where the vector $\bZ_n$ is
%
\begin{eqnarray}\label{zr}
\bZ_{n,r}& =& \sum_{j=1}^{n_r}
\bq(x_{r j}) - \sum_{k, j}
h_r(x_{kj}) \bq(x_{kj} )
\nonumber
\\[-8pt]
\\[-8pt]
\nonumber
&= &\sum
_{k, j} \bigl[ \bigl\{ \delta_{kr} -
h_r(x_{kj}) \bigr\} \bq(x_{kj} ) \bigr],
\end{eqnarray}
and $n^{-1/2}\bZ_n$
is asymptotic normal with mean $\0$ and a positive definite
covariance matrix.
\end{lemma}

Note that $\bZ_n$ is a sum of independent random
vectors with finite moments. The mean of each is not zero,
but the total is zero. In Theorem~\ref{asyofdistest}
we defined the $r$th segment of $\bB_k= \bB_k(\infty)$ as
\begin{eqnarray*}
\bB_{k, r} &=& {\rho_k} \sE\bigl\{\delta_{kr}
\bq(x_{kj}) - h_r(x_{kj}) \bq(x_{kj} )
\bigr\}
\\
&=& \int\bigl\{ \delta_{kr}h_k(x)
- h_r(x ) h_k(x)\bigr\} \bq(x) \,d\bar G(x).
\end{eqnarray*}
An interesting and useful
observation is that { for $1\leq k, r\leq m$,
$\bB_{k, r} =\bB_{r, k}$.}
From $\sum_{k=0}^m h_k(x) = 1$, we may also verify that
\[
\sE\{Z_{n, r}\} = n \sum_{k=0}^m
\bB_{k,r} = 0.
\]

Thus, $n^{-1/2}\bZ_n$ is asymptotic normal
with mean $\0$ and some variance matrix.
This fact together with the form of
the quadratic approximation implies
that $\ell_n(\btheta) - \ell_n(\btheta^*)$ is maximized at the $\hat
{\btheta}$
that satisfies
\[
n^{1/2} \bigl( \hat{\btheta} - \btheta^* \bigr) = n^{-1/2}
\bW^{-1}\bZ_n +o_p(1).
\]
See \citet{HjoPol93} for this justification.

The remaining task is to verify that the asymptotic variance
of $n^{-1/2} \bZ_n$ is given by $\bW-\bW\bS\bW$.
This proves Theorem~\ref{thetaest}.

\subsection{\texorpdfstring{Proof of Theorem \protect\ref{bahadur}}{Proof of Theorem 3.1}}\label{appa2}

The key to the proof is to show a seemingly obvious claim:
$
\hat\xi_{r} - \xi_{r} = O_p
(n^{-1/2} ).
$
This is an immediate consequence of
%
\begin{equation}
\label{G-order} \hat G_r(x) - G_r(x) =
O_p\bigl(n^{-1/2}\bigr)
\end{equation}
uniformly for $x$ in a neighborhood of $\xi_r$.
We now prove (\ref{G-order}). Recall that
\[
\hat G_{r}( x) = n_r^{-1} \sum
_{k,j} h_r(x_{kj}; \hat{
\btheta})I(x_{kj} \leq x).
\]
Replacing $\hat{\btheta}$ in $\hat G_r(x)$
by its true value $\btheta^*$, we define
\[
\tilde G_{r}( x) = n_r^{-1} \sum
_{k,j} { h_r \bigl(x_{kj}; \btheta^*
\bigr)} I(x_{kj}\leq x),
\]
a sum of independent random variables.
From $\sE\{ \tilde G_{r}(x)\} = G_r(x)$, we get
\[
\sE\bigl\{ \tilde G_r(x) - G_{r}(x) \bigr
\}^2 = n_r^{-2} \sum
_{k, j} \var\bigl\{ { h_r\bigl(x_{kj};
\btheta^* \bigr)} I(x_{kj}\leq x) \bigr\}.
\]
Since $0 \leq h_r(x; \btheta^* ) \leq1$,
we have $\var\{ h_r(x_{kj}; \btheta^* ) I(x_{kj} \leq x) \} \leq1$.
Hence,
\[
\sup_x \sE\bigl\{ G_r(x) - \tilde
G_{r}(x) \bigr\}^2 \leq n/ n_r^{2}
= O\bigl(n^{-1}\bigr).
\]
Since $G_r(x)$ and $\tilde G_{r}(x)$ are distribution
functions, the above rate is uniform in $x$.
Hence
\[
\sup_x \bigl| \tilde G_r(x) - G_{r}(x)\bigr|
= O_p \bigl(n^{-1/2} \bigr).
\]
Therefore, (\ref{G-order}) is implied by
$\sup_x |\hat G_r(x) - \tilde G_{r}( x) | = O_p(n^{-1/2})$.
Note that
\begin{eqnarray*}
\tilde G_r(x) - \hat G_{r}(x) &=& n_r^{-1}
\sum_{k,j} \bigl\{ {h_r
\bigl(x_{kj}; \btheta^* \bigr)} - h_r(x_{kj};
\hat{\btheta}) \bigr\} I(x_{kj}\leq x).
\end{eqnarray*}
The partial derivative of $h_r(x; \btheta)$ with
respect to $\btheta$ is bounded by $ \sqrt{m} \|\bq(x)\|$. Thus
\[
\bigl| \tilde G_r(x) - \hat G_{r}(x)\bigr| \leq(n
\rho_r)^{-1} \bigl\|\hat{\btheta} - \btheta^*\bigr\| \sum
_{k,j} \sqrt{m}\bigl\|\bq(x_{kj}) \bigr\|.
\]
The conditions of Theorem~\ref{thetaest} imply that
$\bq(x_{kj})$ has finite moments of any order. Thus,
$
n^{-1} \sum_{k,j} \|\bq(x_{kj})
\| = O_p(1),
$
and subsequently,
\begin{eqnarray*}
\bigl| \tilde G_r(x) - \hat G_{r}(x)\bigr|& \leq&(n
\rho_r)^{-1} \bigl\|\hat{\btheta} - \btheta^*\bigr\| \sum
_{k,j} \sqrt{m}\bigl\|\bq(x_{kj}) \bigr\| = O_p
\bigl( \bigl\|\hat{\btheta} -\btheta^*\bigr\|\bigr) \\
&=& O_p\bigl(n^{-1/2}
\bigr).
\end{eqnarray*}
This completes the proof of (\ref{G-order}).

The classical Bahadur representation
was a rate result in the mode of ``almost sure.'' Our result is stated
in terms
of ``in probability,'' and therefore it has a simpler proof.
As for the classical case, the
representation is equivalent to the following lemma:

%
\begin{lemma}
\label{baha}
Under the conditions of Theorem~\ref{bahadur},
for any $c > 0$ and $r=0, 1, \ldots, m$, we have
\[
\sup_{x\dvtx|x-\xi_r|< cn^{-1/2}}\bigl |\bigl\{ \hat G_r(x) - \hat
G_r(\xi_r)\bigr\} - \bigl\{ G_r(x) -
G_r(\xi_r)\bigr\}\bigr | = O_p
\bigl(n^{-3/4} \bigl(\log(n)\bigr)^{1/2}\bigr).
\]
\end{lemma}

\begin{pf}
We prove this lemma for $r=0$; the other cases are equivalent.
Without loss of generality we assume $x \geq\xi_r$.
Note that
\begin{eqnarray*}
&&\bigl\{\hat G_0(x) - \hat G_0(\xi_r)
\bigr\} - \bigl\{\tilde G_0(x) - \tilde G_0(
\xi_r) \bigr\} \\
&&\qquad= n_0^{-1} \sum
_{k,j} \bigl\{ h_0(x_{kj}; \hat{
\btheta}) - h_0 \bigl(x_{kj}; \btheta^* \bigr) \bigr\} I(
\xi_r< x_{kj} \leq x).
\end{eqnarray*}
By the mean value theorem and the specific form of $h_0(x; \btheta)$,
we have
\[
\bigl| h_0(x_{kj}; \hat{\btheta}) - h_0
\bigl(x_{kj}; \btheta^* \bigr)\bigr | \leq\sqrt{m}\bigl\|\bq(x_{kj})\bigr\|
\bigl\|\hat{\btheta} - \btheta^*\bigr\|.
\]
From $\sE\{ \|\bq(x_{kj})\| \} < \infty$, we get
$\sE\{ \bq(x_{kj}) I(\xi_r < x_{kj} \leq x) \} = O(n^{-1/2})$
and
\[
\bigl\{\hat G_0(x) - \hat G_0(\xi_r)
\bigr\} - \bigl\{\tilde G_0(x) - \tilde G_0(
\xi_r) \bigr\} = O_p \bigl(n^{-1} \bigr).
\]
With this result, Lemma~\ref{baha} is proved if we show that
\[
\sup_{x\dvtx|x-\xi|< cn^{-1/2}} \bigl| \bigl[\tilde G_0(x) - \tilde
G_0(\xi_r) \bigr] - \bigl[ G_0(x) -
G_0(\xi_r) \bigr]\bigr| = O_p
\bigl(n^{-3/4} \bigl(\log(n) \bigr)^{1/2} \bigr).
\]
Since $\tilde G_0(x)$ is a sum of bounded random variables and
$\sE\{\tilde G_0(x)\} = G_0(x)$, the result can be proved
following Lemma 2.5.4E in \citet{Ser80}, page 97;
we omit the details here.
This completes the proof.
\end{pf}

\begin{pf*}{Proof of Theorem \protect\ref{bahadur}}
We have $\hat\xi_{r}-\xi_{r} = O_p(n^{-1/2})$ for any $\xi_r$,
and
the derivative of $G_r$ is positive and continuous in a neighborhood of
$\xi_r$. Therefore,
\[
G_r(\hat\xi_{r}) - G_r(\xi_{r})
= g_r(\xi_{r}) (\hat\xi_{r}-
\xi_{r}) + O_p\bigl(n^{-1}\bigr).
\]
By definition, we have $\hat G_r(\hat\xi_{r}) = \alpha_r + O(n^{-1})$.
Thus, replacing $x$ by $\hat\xi_r$, and $\xi$ by $\xi_r$,
the result of Lemma~\ref{baha} becomes
\[
\bigl| \bigl\{ \alpha- \hat G_r(\xi_{r})\bigr\} -
g_r(\xi_{r}) (\hat\xi_{r}-\xi_{r})
\bigr| = O_p\bigl(n^{-3/4} \bigl(\log(n)\bigr)^{1/2}
\bigr).
\]
This is equivalent to the conclusion of the theorem.
\end{pf*}

\subsection{\texorpdfstring{Proof of Theorem \protect\ref{asyofdistest}}{Proof of Theorem 3.2}}\label{appa3}
Theorem~\ref{asyofdistest} characterizes the asymptotic joint normality
of a number of
MELE distribution estimates. It is proved by approximating $\hat
G_{r}(x)$ by
a summation of independent random variables.

By the proof of Theorem~\ref{thetaest},
$
\hat{\btheta} - \btheta^* = n^{-1} \bW^{-1} \bZ_n
+ o_p(1).
$
Hence,
\begin{eqnarray*}
\hat G_{r}(x) &=& n_r^{-1} \sum
_{k,j} h_r(x_{kj}; \hat{\btheta})
I(x_{kj}\leq x)
\\
&=& n_r^{-1} \sum_{k,j} \bigl[
h_r(x_{kj}) + \bigl\{\dot{h}_r
\bigl(x_{kj}; \btheta^*\bigr)\bigr\}^\tau\bigl(\hat{\btheta}
- \btheta^*\bigr)\bigr] I(x_{kj}\leq x) + o_p
\bigl(n^{-{1}/{2}}\bigr)
\\
&=& n_r^{-1} \sum_{k,j}
h_r(x_{kj})I(x_{kj}\leq x)
\\
&&{} + n_r^{-1} \biggl\{ n^{-1} \sum
_{k,j} \dot{h}_r \bigl(x_{kj}; \btheta^*
\bigr) I(x_{kj} < x) \biggr\}^\tau\bW^{-1}
\bZ_n + o_p\bigl(n^{-{1}/{2}}\bigr),
\end{eqnarray*}
where
$
\dot{h}_r (x; \btheta^* ) = {\partial h_r
(x; \btheta) }/{\partial\btheta} |_{\btheta= \btheta^*}
$.
Working out the expression of $ \dot{h}_r (x; \btheta^*)$ in terms of
$\bq(x)$ and $h_r(x)$,
and by the law of large numbers, we find that, almost surely,
\[
n^{-1} \sum_{k,j} \dot{h}_r
\bigl(x_{kj}; \btheta^* \bigr) I(x_{kj} < x) \to
\bB_r(x),
\]
where { $\bB_r(x)$ is defined in the theorem. }
We remark here that $ \lim_{x \to\infty} \bB_{r}(x) = \bB_{r}$;
the latter was defined in the proof of Theorem~\ref{thetaest}.
Before the final step, we may verify that
\[
G_r(x) = n_r^{-1} { \sE} \biggl\{ \sum
_{k,j} h_r(x_{kj})I(x_{kj}
\leq x) \biggr\}.
\]
These preparations enable us to write
\begin{eqnarray*}
\hat G_{r}(x) - G_r(x) &=& n_r^{-1}
\sum_{k,j} \bigl[ h_r(x_{kj})I(x_{kj}
\leq x) - \sE\bigl\{ h_r(x_{kj})I(x_{kj}\leq x)
\bigr\} \bigr]
\\
&& {}+ n_r^{-1} \bB_r^{\tau}(x)
\bW^{-1} \bZ_n + o_p\bigl(n^{1/2}
\bigr).
\end{eqnarray*}
The two leading terms are summations of { independent random vectors
and both have mean zero}.
The joint asymptotic normality of
$\sqrt{n}\{\hat G_{r}(x) - G_r(x) \}$
and $\sqrt{n}\{\hat G_{s}(y) - G_s(y) \}$ is hence implied.
We derive the algebraic expression of
$ \omega_{rs}(x,y)$ in the next subsection.

\subsubsection{\texorpdfstring{Asymptotic covariance $\sigma_{rs}(x,y)$}
{Asymptotic covariance sigma rs(x,y)}}
\label{secasycov}
From the expansion of $\hat G_{r}(x) - G_r(x)$,
$ \omega_{rs}(x,y)$ is decomposed into four covariances.
Using $\var(\bZ_n) =n( \bW- \bW\bS\bW)$ as shown earlier,
we find that one of them is given by
\[
\cov\bigl(\bB_r^{\tau}(x)\bW^{-1}
\bZ_n, \bB_s^{\tau}(y)\bW^{-1}
\bZ_n \bigr) = n \bB_r^{\tau}(x) \bigl(
\bW^{-1} - \bS\bigr) \bB_s(y).
\]
We build another term from the following computations:
\begin{eqnarray*}
&& \cov\bigl\{ h_r(x_{kj})I(x_{kj}\leq x),
h_s(x_{kj})I(x_{kj}\leq y)\bigr\}
\\
&&\qquad= \rho_k^{-1} \int_{-\infty}^{(x \wedge y)}
h_r(x) h_s(x) h_k(x) \,d\bar G(x) -
\rho_k^{-2} c_{rk}(x) c_{sk}(y),
\end{eqnarray*}
where
$
c_{rs}(x)
=
\int_{-\infty}^x h_r(t) h_s(t) \,d\bar G(t)$.
Because $\sum_{k=0}^m h_k(x) = 1$, we get
\begin{eqnarray*}
&& \sum_{k=0}^m n_k \cov
\bigl\{ h_r(x_{kj})I(x_{kj}\leq x),
h_s(x_{kj})I(x_{kj}\leq y)\bigr\}
\\
&&\qquad= n\Biggl\{ c_{rs}(x\wedge y) - \sum_{k=0}^m
\rho_k^{-1} c_{rk}(x) c_{sk}(y)\Biggr
\}.
\end{eqnarray*}
The last task is the cross-term $\cov\{ h_r(x_{kj})I(x_{kj}\leq x), \bZ
_n\}$.
We break $\bZ_n$ into segments $\bZ_{n, s}$ and
then into centralized $\tilde{\bZ}_{n, s}$.
\begin{eqnarray*}
&&\cov\bigl\{ h_r(x_{kj})I(x_{kj}\leq x),
\bZ_{n, s} \bigr\}\\
&&\qquad= \sE\bigl[ h_r(x_{kj})I(x_{kj}
\leq x) \bigl\{ { \bigl(\delta_{kr} - h_s(x_{kj})
\bigr)}\bq^\tau(x_{kj}) - \rho_k^{-1}
\bB_{k, s}^\tau\bigr\} \bigr]
\\
&&\qquad= \rho_k^{-1} { \bB^{\tau}_{k, r}(x)} -
\rho_k^{-2} c_{rk}(x) { \bB^{\tau}_{k,s}}.
\end{eqnarray*}
Summing over $\{k, j\}$, the first term sums to zero, so we find
\[
\cov\biggl\{ \sum_{k, j} h_r(x_{kj})I(x_{kj}
\leq x), \bZ_{n, s} \biggr\} = - n \sum_{k=0}^m
\rho_k^{-1} c_{rk}(x) {\bB^{\tau}_{k,s}}.
\]
Next, we assemble $\bB_{k, s}$ over $s$ to get $\bB_k$,
\[
\cov\biggl\{ \sum_{k, j} h_r(x_{kj})I(x_{kj}
\leq x), \bZ_{n} \biggr\} = - n \sum_{k=0}^m
\rho_k^{-1} c_{rk}(x) { \bB^{\tau}_{k}
}.
\]
Entering $ \bB_{s}^\tau(y) \bW^{-1}$ into the second argument of the covariance,
we get
\[
\cov\biggl\{ \sum_{k, j} h_r(x_{kj})I(x_{kj}
\leq x), { \bB_{s}^\tau(y)} \bW^{-1}
\bZ_{n}\! \biggr\} = - n \sum_{k=0}^m
\rho_k^{-1} c_{rk}(x) \bB_{k}^\tau
\bW^{-1} \bB_{s}(y).
\]
Thus, the covariance between $\sqrt{n}\{\hat G_{r}(x) - G_r(x) \}$
and $\sqrt{n}\{\hat G_{s}(y) - G_s(y) \}$ is given by
\begin{eqnarray*}
\omega_{rs}(x,y) &= & (\rho_r\rho_s)^{-1}
\Biggl[ c_{rs}(x \wedge y) - \sum_{k=0}^m
\rho_k^{-1} c_{rk}(x) c_{sk}(y)
\\
&&\hspace*{43pt}{} - \sum_{k=0}^m \rho_k^{-1}
\bB_k^\tau\bW^{-1} \bigl\{ { c_{rk}(x)
\bB_s(y) + c_{sk}(y) \bB_r(x)} \bigr\}\\
&&\hspace*{117pt}\qquad{} +
\bB_r^{\tau}(x) \bigl(\bW^{-1} - \bS\bigr)
\bB_s(y) \Biggr].
\end{eqnarray*}

Further simplification is possible.
We find
$c_{r0}(x) + \bB_r^{\tau}(x)\bW^{-1} \bB_0 = 0$,
$c_{rk}(x) + \bB_r^{\tau}(x)\bW^{-1} \bB_k = \delta_{rk} G_r(x)$
and
$ \sum_{k=0}^m \rho_k^{-1} \bW^{-1} \bB_k\bB_k^{\tau} \bW^{-1} = \bS$.
These findings lead to
\[
\omega_{rs}(x,y) = (\rho_r\rho_s)^{-1}
\bigl\{ c_{rs}(x \wedge y) + \bB_r^{\tau}(x)
\bW^{-1}\bB_s(y) - \rho_r \delta_{rs}
G_r(x)G_s(y) \bigr\}.
\]
The final expression (\ref{omegars}) is obtained by
noticing that
$
c_{rs}(x)
=
\rho_r \delta_{rs} G_r(x) - a_{rs}(x),
$
where $a_{rs}(x)$ was defined in Theorem~\ref{asyofdistest}.

\subsection{\texorpdfstring{Proof of Theorem \protect\ref{matdiff}}{Proof of Theorem 4.1}}\label{appa4}
Both the EL and EM quantiles admit Bahadur representations, and
it suffices to show the same conclusion for the distribution estimators.
For the bivariate case, we denote the
asymptotic covariance matrices of the EL and EM distributions as
\begin{eqnarray*}
\Omega_{\EL} = \pmatrix{ \sigma_{rr}(x,x)
& \sigma_{rs}(x,y)
\vspace*{2pt}\cr
\sigma_{rs}(x,y) & \sigma_{ss}(y,y)
},\qquad \Omega_{\EM} = \pmatrix{
\omega_{rr}(x,x) & \omega_{rs}(x,y)
\vspace*{2pt}\cr
\omega_{rs}(x,y) & \omega_{ss}(y,y)},
\end{eqnarray*}
where $x=\xi_r $ and $y=\xi_s$ are two population quantiles or two real values.
We show that $\Omega_{\EM} - \Omega_{\EL} $
is nonnegative definite by writing it as
$D_{11} - D_{12}D_{22}^{-1}D_{21}$, with the $D_{ij}$ being blocks
of a nonnegative definite matrix $D$.
By standard matrix theory, the nonnegative definiteness of $D$ implies
that of $D_{11} - D_{12}D_{22}^{-1}D_{21}$.
The generic element of $\Omega_{\EM} - \Omega_{\EL} $ is
$
a_{rs}(x, y) - \bB_r^{\tau}(x) \bW^{-1}
\bB_s(y),
$
which fits into
$D_{11} - D_{12}D_{22}^{-1}D_{21}$
with
\[
D_{11} = \pmatrix{ a_{rr}(x) &
a_{rs}(x \wedge y)
\vspace*{2pt}\cr
a_{rs}(x \wedge y) & a_{ss}(y)},\qquad
D_{12} = \pmatrix{ \bB_r^{\tau}(x)
\vspace*{2pt}\cr
\bB_s^{\tau}(y)},
\]
$D_{21} = D_{12}^{\tau}$ and $D_{22} = \bW$.
We will show that $D = \int\mathbf{U}(z) \,d \bar G(z)$
for some nonnegative definite $\mathbf{U}(z)$ for all $z$.
Then $D$ is nonnegative definite and so is
$
\Omega_{\EM} - \Omega_{\EL}
=
D_{11} - D_{12}D_{22}^{-1}D_{21}.
$

We now search for such a $\mathbf{U}(z)$.
We write
\begin{eqnarray*}
a_{rs} (x \wedge y ) &=& \int I( z \leq x) I( z \leq y)\bH[
r+1, s+1](z) \,d\bar G(z),
\\
\bB_{r} (x) &=& \int I( z \leq x) \bH[-1, r+1](z)
\otimes\bq(z) \,d\bar G(z),
\\
\bW&=& \int\bH[-1, -1](z) \otimes\bigl\{ \bq(z)
\bq^{\tau}(z) \bigr\} \,d\bar G(z).
\end{eqnarray*}
Using the Khatrin--Rao product operator $\ast$ [\citet{LiuTre08}],
we find such a $\mathbf{U}(z)=A_1(z) \ast A_2(z)$ with
\begin{eqnarray*}
A_1(z) &= & P \pmatrix{
 \bH[r+1,
r+1](z) & \bH[r+1, s+1](z) & \bH[r+1,
-1](z)
\vspace*{2pt}\cr
\bH[s+1, r+1](z) & \bH[s+1, s+1](z)
& \bH[
s+1, -1](z)
\vspace*{2pt}\cr
\bH[-1, r+1](z) & \bH[-1, s+1](z) & \bH[
-1, -1](z)} P,
\\
A_2(z) &=& \pmatrix{ 1 & 1 &
\bq^{\tau}(z)
\vspace*{2pt}\cr
1 & 1 & \bq^{\tau}(z)
\vspace*{2pt}\cr
\bq(z) & \bq(z) & \bq(z)\bq^{\tau}(z)},\qquad  P=
\pmatrix{ I( z \leq x ) & 0 & 0
\vspace*{2pt}\cr
0 & I( z \leq y ) & 0
\vspace*{2pt}\cr
0 & 0 & \bI_m
}.
\end{eqnarray*}
The matrix $A_2(z)$ is clearly nonnegative definite for any $z$.
Note that $\bH(z)$ is nonnegative definite for any $z$;
the nonnegative definiteness of $A_1(z)$ is an easy consequence.
Since the $\ast$ product of two nonnegative definite matrices
is still nonnegative definite [Lemma 5 of \citet{LiuTre08}],
we conclude that $\mathbf{U}(z) = A_1(z) \ast A_2(z)$
is also nonnegative definite for any $z$.
This completes the proof. This proof can easily be extended to the case where
more distributions or quantiles are involved.
\end{appendix}

\section*{Acknowledgments}
We are grateful to the referees, the Associate Editor, and the Editor
for helpful comments.



%

\printaddresses

\end{document}